\documentclass{amsart}
\usepackage{graphicx}
\newtheorem{thm}{Theorem}[section]
\newtheorem{cor}[thm]{Corollary}
\newtheorem{lem}[thm]{Lemma}
\newtheorem{prop}[thm]{Proposition}
\theoremstyle{definition}
\newtheorem{defn}[thm]{Definition}
\theoremstyle{remark}
\newtheorem{rem}[thm]{Remark}
\numberwithin{equation}{section}

\newcommand{\RR}{\mathbb R}
\newcommand{\CC}{\mathbb C}

\newcommand{\C}{\mathcal{C}}

\newcommand{\dc}{d^c}

\DeclareMathOperator{\RE}{Re}
\DeclareMathOperator{\IM}{Im}
\DeclareMathOperator{\trace}{trace}

\DeclareMathOperator{\Hess}{Hess}

\newcommand{\zz}[1]{\frac{\partial^2#1}{\partial z\partial {\bar z}}}
\newcommand{\ww}[1]{\frac{\partial^2#1}{\partial w\partial {\bar w}}}
\newcommand{\zw}[1]{\frac{\partial^2#1}{\partial z\partial {\bar w}}}

\begin{document}

\title{On Stein Neighborhood Basis of Real Surfaces}
\author{Marko Slapar}%
\address{Institute of Mathematics, Physics and Mechanics, University of Ljubljana, Jadranska 19, 1000 Ljubljana, Slovenia}%
\email{marko.slapar@fmf.uni-lj.si}%
\date{September 2003}%
\begin{abstract}
In this paper, we show that a compact real surface
embedded in a complex surface has a regular Stein neighborhood
basis, provided that there are only finitely many complex points
on the surface, and that they are all flat and hyperbolic.
An application to unions of totally real planes in $\CC^2$ is then given.
\end{abstract}
\maketitle
\section{Introduction}
Let $S\hookrightarrow X$ be a real compact surface, smoothly
embedded in a complex surface $X$. In a generic position, there
are only finitely many complex points on $S$, which can be,
following Bishop \cite{Bis}, classified as either elliptic or
hyperbolic. \par By a result of Bishop \cite{Bis}, the local
holomorphic hull of $S$ at elliptic points contains a real one
parameter family of pairwise disjoint holomorphic discs with
boundaries in $S$. The presence of elliptic complex points on $S$
thus prevents the surface $S$ from having small Stein
neighborhoods in $X$. On the other hand, by a result of
Forstneri\v c and Stout \cite{FoS}, $S$ is locally polynomially
convex at hyperbolic complex points. Using this, one can easily
construct a Stein neighborhood basis of a real surface with only
hyperbolic complex points by just patching together local
pseudoconvex defining functions. The problem is that one does not
necessarily understand the topology of such a basis. For example,
none of the members are a priori even of the same homotopy type as
$S$. For this reason, we prefer to put some further restrictions
on our basis.
\begin{defn} Let $\pi\!: S\hookrightarrow M$ be an embedding
of a manifold $S$ into a manifold $M$. A system
$\{\Omega_\epsilon\,,\epsilon\in(0,1)\}$ of open neighborhoods of
$\pi(S)$ in $M$ is called a \emph{regular basis}, if for every
$\epsilon\in (0,1)$, we have
\begin{enumerate}
\item$\Omega_\epsilon=\bigcup_{s< \epsilon}\Omega_s$, \item
$\bar\Omega_\epsilon=\bigcap_{s>\epsilon}\Omega_s$, \item
$\pi(S)=\bigcap_{s>0}\Omega_s$ is a strong deformation retract of
$\Omega_\epsilon$.
\end{enumerate}
\end{defn}
It has been asked by Forstneri\v c in \cite{For2}, whether the
nonexistence of elliptic complex points on $S$ is sufficient for
construction a \emph{regular} Stein neighborhood basis of $S$ in
$X$. We partially answer the question in the following theorem.
\begin{thm}\label{maintheorem} Let $S$ be a compact real surface, embedded in a complex
surface $X$, and having only finitely many complex points. Let all
complex points on $S$ be flat hyperbolic complex points. Then $S$
has a regular strictly pseudoconvex Stein neighborhood basis in
$X$.
\end{thm}
The condition of flatness, with the rest of the terms, will be
defined in the next section. We would like to indicate that the
flatness condition is in our view redundant, and is here only
because of our inability to yet prove the theorem without it.
\section{Local structure at complex points}
Let $S\hookrightarrow X$ be a real compact surface, embedded in a
complex surface $X$. We call $p\in S$ a \emph{complex} point if
$T_pS\subset T_pX$ is a complex line. If this is not the case, the
point $p$ is \emph{totally real}. In a generic situation, there
are only finitely many complex points on $S$. Let $p\in S$ be an
isolated complex point on $S\hookrightarrow X$. In local
holomorphic coordinates around $p=(0,0)$, we can write $S$ as a
graph $w=f(z)$. In a generic case, we can further change the
coordinates, see \cite{Bis}, to locally have $S$ given by
\begin{equation}\label{eqn1}
2w=\alpha z\bar z +\frac{1}{2}(z^2+{\bar z}^2)+o(|z|^3),
\end{equation}
where $0\le\alpha\le \infty$. The case $\alpha=\infty$ should be
understood as the surface $w=z\bar z+o(|z|^3).$ We call $p$ a
\emph{hyperbolic} complex point if $0\le\alpha<1$ and
\emph{elliptic} if $\alpha>1$. The parabolic case $\alpha=1$ is
not generic. A complex point $p$ on a real surface $S$ is called
\emph{quadratic}, if we can write the surface in the above form
with $o(|z|^3)\equiv 0$. We call a complex point $p$ a
\textit{flat} complex point, if local holomorphic coordinates can
be chosen so that $\IM o(|z|^3)\equiv 0$ in the form (\ref{eqn1}).
Quadratic complex points are of course flat.
\begin{rem} If $S$ is embedded as a real analytic submanifold of a complex surface $X$,
then by a result of Moser and Webster in \cite{MoW}, every real
analytic elliptic complex point is flat. On the other hand, they
show the existence of real analytic hyperbolic complex points that
are not flat.
\end{rem}
Let us now take the case of a flat hyperbolic complex point.
Locally, we assume the surface lies in $\CC^2$ with coordinates
$(z,w)$, the complex point corresponds to $(z,w)=(0,0)$ and that
the surface $S$ is given by
\begin{equation}\label{formatocke}
S=\left\{ \begin{array}{ll} \RE w &= \frac{1}{2}\alpha z\bar z +
\frac{1}{4}(z^2+{\bar z}^2)+\tau\\
\IM w &=0\end{array}\right.
\end{equation}
where $\tau=o(|z|^3)$ is a real function, and $0\le\alpha<1$. We
introduce new real coordinates in the neighborhood of the origin
of $\CC^2$, given by
\begin{equation}\label{nonholomorphic}\begin{array}{ll}x&=\RE z\\y&=\IM z\\u&=\RE
w-\frac{1}{2}\alpha z\bar z-\frac{1}{4}(z^2+{\bar
z}^2)-\tau\\v&=\IM w.\end{array}\end{equation} These new
coordinates are nonholomorphic and they depend on $\alpha$, with
the surface $S$ corresponding to $\{u=v=0\}$. The partial
derivatives in the Levi form computation are expressed in these
coordinates as
\begin{equation}\label{changelevi}\begin{array}{ll}
4\frac{\partial^2}{\partial z\partial\bar z}&=\triangle_{x,y}
-2((\alpha+1)x\frac{\partial}{\partial
x}+(\alpha-1)y\frac{\partial}{\partial y} +
\alpha)\frac{\partial}{\partial u}
\\&+((\alpha+1)^2x^2+(\alpha-1)^2y^2)\frac{\partial^2}{\partial
u^2}\\
&+ (o(|z|^2)\frac{\partial }{\partial x}+o(|z|^2)\frac{\partial
}{\partial y}+o(|z|))\frac{\partial }{\partial
u}+o(|z|^3)\frac{\partial^2
}{\partial u^2}\\
4\frac{\partial^2}{\partial w\partial\bar w}=&\triangle_{u,v}\\
4\frac{\partial^2}{\partial z\partial\bar w}=&
\left(\frac{\partial^2}{\partial x\partial
u}-(\alpha+1)x\frac{\partial^2 }{\partial
u^2}+o(|z|^2)\frac{\partial^2}{\partial
u^2}\right)\\
&+i\left(-\frac{\partial^2}{\partial y\partial
u}+(\alpha-1)y\frac{\partial^2}{\partial u^2}
+o(|z|^2)\frac{\partial^2}{\partial u^2}\right).
\end{array}
\end{equation}
The $o(|z|), o(|z|^2)$ and $o(|z|^3)$ terms consist of derivatives
of $\tau$, and so vanish for quadratic complex points. The term
$\triangle_{x,y}$ denotes the Laplace operator with respect to
$x,y$ coordinates.
\section{Local construction at hyperbolic complex points}
Throughout this section, let $(z,w)$ be the standard holomorphic
coordinates coordinates on $\CC^2$ and let $(x,y,u,v)$ be the
non-holomorphic real coordinates on $\CC^2$, given by
(\ref{nonholomorphic}), and are dependent on $\alpha$. \par In the
construction of Stein neighborhood basis near a hyperbolic complex
point, we use the following lemma.
\begin{lem}\label {mama lemma} Let $S$ be a submanifold in $M$ and let
$\Omega$ be an open neighborhood of $S$ in $M$. Let
$f\!:\Omega\rightarrow[0,1)$ be a $\C^2$ function, having the
property $S=\{f=0\}=\{\nabla f=0\}$. Then
$\{\Omega_\epsilon\}_{\epsilon<1}=\{f<\epsilon\}_{\epsilon<1}$,
defines a regular neighborhood basis for $S$ in $M$.
\end{lem}
\begin{proof} Let $\Omega_s=\{f< s\}$, and let $\psi_t$ be a the
flow of the vector field $-\nabla f$ in some Riemannian metric on
$M$. This flow gives us a strong deformation retract of $\Omega_s$
to $S$.
\end{proof}
\par Our goal is to
find functions $\phi(x,y,u)\ge 0$ and $\psi(x,y,v)\ge 0$, defined
in a small neighborhood of $(0,0,0)$ in $\RR^3$ and having the
following properties
\begin{itemize} \item $\Phi(z,w)=\phi(x,y,u)+\psi(x,y,v)$ is
plurisubharmonic in a small neighborhood $U$ of $(0,0)$ in
$\CC^2$, and strictly plurisubharmonic in $U\backslash\{(0,0)\}$,
\item $\{\Phi=0\}=\{\nabla \Phi=0\}=\{u=v=0\}\bigcap U$.
\end{itemize}
\subsection{The quadratic case}
We first study quadratic complex points. Let $S$ be written as in
(\ref{formatocke}) with vanishing $\tau$. We are frequently going
to use the following simple observation
\begin{lem} \label{pos}
Let $p(x,y,u)=b_2(x,y)u^2+b_1(x,y)u+b_0$, where $b_0,b_1,b_2$ are
continuous functions, defined in a neighborhood of the origin.
Assume $b_1,b_0$ both vanish at $(x,y)=(0,0)$ and $b_2>0$. Then
there exists a small neighborhood $U$ of $(0,0,0)\in \RR^3$ with
$p$ strictly positive on $U\backslash\{u=0\}$, as long as
$b_1^2<4b_2b_0$ for small $(x,y)\ne (0,0)$.
\end{lem}
\begin{proof} For fixed $x,y$,
$$\frac{-b_1\pm\sqrt{b_1^2-4b_2b_0}}{2b_2}$$
gives the zeros of the quadratic polynomial
$b_2(x,y)u^2+b_1(x,y)u+b_0$. We notice that when $(x,y)$
approaches $(0,0)$, both zeros converge to $0$. The only way to
achieve strict positivity of $p$ in some $U\backslash \{u=0\}$, is
for both zeros to be in $(\CC\backslash \RR)\cup \{0\}$ for small
$(x,y)$. That is true if $b_1^2<4b_2b_0$ for small nonzero
$(x,y)$.
\end{proof}
\begin{lem}\label{L1} For $\alpha < 0.52$, and any $M>0$,
there exists a homogeneous polynomial $P\in\RR[x^2,y^2,u]$ of
degree $6$ and an open neighborhood $U\subset \CC^2$ of $(0,0)$,
so that the function $$
\Phi(z,w)=P(x^2,y^2,u)+M(x^2+y^2)u^6+(1+x^2+y^2)v^2$$ has
properties
\begin{enumerate}
\item[(a)]$\Phi$ is strictly plurisubharmonic for $(z,w)\in
U\backslash\{(0,0)\}$, \item[(b)]
$(u,v)\cdot(\frac{\partial\Phi}{\partial u},\frac{\partial
\Phi}{\partial v})>0$ for $(z,w)\in U\backslash\{u=v=0\}$,
 \item[(c)] $\Phi>0$ for $(z,w)\in U\backslash\{u=v=0\}$,
\item[(d)] $\{\Phi=0\}=\{u=v=0\}\cap U$.
\end{enumerate}
\end{lem}
\begin{proof}
We look at polynomials ${\tilde P}$ of the form
$$\begin{array}{ll} {\tilde P}(x^2,y^2,u)&=u^6+((6\alpha+c)x^2-cy^2)u^5+(Ax^4+Bx^2y^2+A'y^4)u^4\\
&+(Cx^4y^2)u^3+(Dx^6y^2+Ex^4y^4)u^2,\end{array}$$ where
$c,A,A',B,C,D,E$ are strictly positive constants, to be determined
later. The computation of $\frac{\partial^2 {\tilde P}}{\partial
z\partial \bar z}$, using (\ref{changelevi}), gives us
\begin{equation}\label{haha}\begin{array}{ll}4\frac{\partial^2 {\tilde P}}{\partial
z\partial \bar z}&=[(30(\alpha+1)^2+12A+2B-60\alpha(3\alpha+2)a-10c(3\alpha+2))x^2\\
&+(30(\alpha-1)^2+12A'+2B+10c(\alpha-1)^2)y^2]u^4\\
&+[(120\alpha(\alpha+1)^2+20c(\alpha+1)^2+2C-8A(5\alpha+4))x^4\\
&+(120\alpha(\alpha-1)^2-80\alpha c+ 12C-40B\alpha)x^2y^2\\
&+(-20c(\alpha-1)^2-8A'(5\alpha-4))y^4]u^3+[(12A(\alpha+1)^2+2D)x^6\\
&+(12A(\alpha-1)^2+12B(\alpha+1)^2+12A'(\alpha+1)^2+30D+12E-6C(7\alpha+2))x^4y^2\\
&+(12A'(\alpha+1)^2+12B(\alpha-1)^2+12E)x^2y^4+(12A'(\alpha-1)^2)y^6]u^2\\
&+[(6C(\alpha+1)^2-4D(9\alpha+4))x^6y^2\\
&+(6C(\alpha-1)^2-36E\alpha)x^4y^4]u+2(Dx^6y^2+Ex^4y^4)((\alpha+1)^2x^2+(\alpha-1)^2y^2).
\end{array}
\end{equation}
We want this expression to be nonnegative. To insure this, we set
the coefficients so that polynomials in front of $u^3,u^1$ vanish,
and polynomials in front of powers $u^4,u^2,u^2$ are positive for
$(x,y)$ small and nonzero.
\par We first assume that $c=0$, and so also $A'=0$. Setting polynomials in front of
powers $u^3,u^1$ to zero, gives us the following conditions
\begin{enumerate} \item[1.]
$A=15\frac{\alpha(\alpha+1)^2}{5\alpha+4}+\frac{C}{4(5\alpha+4)},$
\item[2.] $B=3(\alpha-1)^2+\frac{3}{10}\frac{C}{\alpha}$,
\item[3.] $D=\frac{3}{2}\frac{(\alpha+1)^2}{9\alpha+4}C$,
\item[4.] $E=\frac{1}{6}\frac{(\alpha-1)^2}{\alpha}C$.
\end{enumerate}
Next we write inequalities, insuring the polynomials in front of
$u^4,u^2,u^0$ are nonnegative. Satisfying the obvious requirement
$C>0$, and using equalities $1-4$, we only need to satisfy the
following inequalities
\begin{enumerate}
\item[5.] $15(\alpha+1)^2+6A+B>30\alpha(3\alpha+2)$, \item[6.]
$2A(\alpha-1)^2 +2B(\alpha+1)^2 +5D+2E>C(7\alpha+2)$.
\end{enumerate}
We also want ${\tilde P}>0$ and $\frac{\partial {\tilde
P}}{\partial u}\ne 0$, when $u\ne 0$. Let us, still assuming
$c=0$, split ${\tilde P}={\tilde P}_1+{\tilde P}_2$ with
$$\begin{array}{l}{\tilde P}_1=u^6+6\alpha x^2u^5+Ax^4u^4\\
{\tilde P}_2=Bx^2y^2u^4+Cx^4y^2u^3+(Dx^6+Ex^4y^4)u^2.\end{array}$$
We require ${\tilde P}_1$ and $u\cdot(\frac{\partial {\tilde
P}_1}{\partial u})$ to be strictly positive when $u\ne 0$, and
${\tilde P}_2$ and $u\cdot(\frac{\partial {\tilde P}_2}{\partial
u})$ nonnegative when $u\ne 0$. Applying Lemma \ref{pos}, we get
two new inequalities
\begin{enumerate}
\item[7.] $75\alpha^2<8A$, \item[8.] $9C^2<32BD$.
\end{enumerate}
We want to find positive constants $A,B,C,D,E$, satisfying
equalities and inequalities $1-8$. Applying equalities $1-4$, we
get constants $A,B,D,E$ expressed in terms of $\alpha$ and $C$.
Inequalities $5-8$ then become
\begin{enumerate}
\item[5.] $30\alpha(3\alpha+2)(5\alpha^2+2\alpha-2)<(5\alpha+2)C$
\item[6.]
$180\alpha(9\alpha+4)(5\alpha+2)(\alpha-1)^2(\alpha+1)^2>
(7\alpha+2)(495\alpha^3+518\alpha^2-64\alpha-112)C$ \item[7.]
$15\alpha(17\alpha^2+4\alpha-8)<2C$ \item[8.]
$80\alpha(\alpha-1)^2(\alpha+1)^2>(37\alpha^2+4\alpha-8)C.$
\end{enumerate}
We can quickly see that, for $\alpha<1$, inequality $7$ follows
from $5$. Let us look at the polynomial
$$q_6=495\alpha^3+518\alpha^2-64\alpha-112$$
from the right hand side of inequality $6.$ We have $q_6(-1)<0,\
q_6(-0.7)>0$ and $q_6(0)<0$, and so the polynomial $q_6$ has three
real zeros, out of which only one, $\alpha_0$, is in the interval
$[0,1].$ Since $q_6(0.43)<0$ and $q_6(0.44)>0$, this zero is in
the interval $(0.43,0.44)$. Applying the quadratic formula, we see
that the polynomial
$$q_8=37\alpha^2+4\alpha-8$$ from the
right hand side of $8$ has two real zeros, of which only one,
$0.41<\alpha_1<0.42$, is positive. For $\alpha<\alpha_1$, both
equations $6$ and $8$ are trivially satisfied with any positive
$C$, since the coefficients in front of $C$ are non-positive. For
$\alpha\le\alpha_0$, we only need to take $C$ small enough to
satisfy $8$, since $5$ is trivially satisfied there for any
positive $C$. This is true because the left hand side of $5$ is
negative for $\alpha<\alpha_0$, since the positive zero
$\frac{-2+\sqrt{24}}{10}$ of $5\alpha^2+2\alpha-2$ is greater than
$\alpha_0$. \par For $0.52>\alpha>\alpha_0$, $6$ follows from $8$.
To see this, we need to check that
$$\frac{180\alpha(9\alpha+4)(5\alpha+2)(\alpha-1)^2(\alpha+1)^2}
{(7\alpha+2)(495\alpha^3+518\alpha^2-64\alpha-112)}>\frac{80\alpha(\alpha-1)^2(\alpha+1)^2}
{(37\alpha^2+4\alpha-8)}$$ holds for $\alpha_0<\alpha<0.52$. This
is equivalent to
$$\begin{array}{l}
9(9\alpha+4)(5\alpha+2)(37\alpha^2+4\alpha+8)-4(7\alpha+2)(495\alpha^2+4\alpha+8)\\
=5(\alpha-4)(225\alpha^3+62\alpha^2-64\alpha-16)>0.\end{array}$$
Let $q=225\alpha^3+62\alpha^2-64\alpha-16$. Then
$q'=675\alpha^2+124\alpha-64$, and $q'$ has only one positive
zero, as we can see from by using the quadratic formula. Since
$q(0)<0$ and the leading coefficient of $q$ is positive, this
implies that $q$ has only one positive zero. Since $q(0.52)<0$,
this zero is greater than $0.52$. So $5(\alpha-4)q>0$ for
$0<\alpha<0.52$, as long as $c$ is small enough.
\par We now only need to worry about $5$ and $8$, for
$\alpha>\alpha_0$. We can simultaneously solve these inequalities
if and only if
$$\begin{array}{ll}
&30\alpha(3\alpha+2)(5\alpha^2+2\alpha-2)(37\alpha^2+4\alpha-8)\\&
\phantom{111111111111}-80\alpha(\alpha-1)^2(\alpha+1)^2(5\alpha+2)<0.\end{array}$$
By estimating the zeros of this polynomial, we see that the
inequality  holds for $\alpha_0<\alpha<0.52$. This shows that
equalities and inequalities $1-8$ can all simultaneously be
satisfied, as long as $0\le \alpha\le 0.52$.
\par By a continuity argument, we can see that condition $c=0$ can
be dropped, as long as we take $c$ small enough and positive. This
forces $A'>0$. The reason is that equalities and inequalities,
insuring positivity of $\frac{\partial^2 {\tilde P}}{\partial
z\partial {\bar z}}$, depend continuously on $c$. Finally, we want
to make ${P}$ more generic in the coefficient in front of $u^2$,
so that $\frac{\partial^2 P}{\partial z\partial {\bar z}}|_{u=0}$
is positive for every $(x,y)\ne (0,0)$. We can simply achieve
this, by taking $P$ to be a small perturbation of $\tilde P$, of
the form $P=\tilde P+\epsilon (x^8+y^8)u^2$. We should again
understand $P$ as a polynomial in variables $x^2$, $y^2$ and $u$.
As long as $\epsilon$ is small enough, we do not spoil the
positivity of $\frac{\partial^2 P}{\partial z\partial {\bar z}}$.
\par Next we compute
$$\frac{\partial^2 (M(x^2+y^2)u^6)}{\partial z\partial {\bar z}}=
4Mu^6-p_1u^5+p_2u^4.$$ Here $p_1,p_2\in \RR[x^2,y^2]$ are
homogeneous polynomials of respectively degrees $1$ and $2$. We
have $\frac{\partial^2 P}{\partial z\partial {\bar
z}}=q_1u^4+q_3u^2+q_5$, where $q_1\in \RR[x^2,y^2]$ is a
polynomial with strictly positive coefficients of degree $1$. For
any small $\epsilon$, the polynomial $$4Mu^6+p_1u^5+(p_2+\epsilon
q_1)u^4$$ is strictly positive for $u\ne 0$ and $(x,y)$ small
enough, and since $$4\frac{\partial ((1+x^2+y^2)v^2)}{\partial
z\partial {\bar z}}= 4v^2\ge 0,$$ we have that $4\frac{\partial
\Phi}{\partial z\partial {\bar z}}> 0$, whenever $(z,w)\ne 0$.
\par From $u\frac{\partial P}{\partial u}>0$ and $P>0$ for $u\ne
0$ we also get $(b)$, $(c)$. This is true, since the term
$M(x^2+y^2)u^6$ is small with respect to $P$.
\par To conclude the proof, we check that $\Phi$ is indeed
strictly plurisubharmonic. First, let us check that in a small
neighborhood of the origin, $\rho=(\frac{1}{2}+x^2+y^2)v^2$ is
plurisubharmonic, and strictly plurisubharmonic if $v\ne 0$. Using
(\ref{changelevi}), we have
$$\begin{array}{ll}
\zz\rho&=4v^2\\
\ww\rho&=2(\frac{1}{2}+x^2+y^2)\\
\zw\rho&=4yv+4ixv.
\end{array}$$
Let us compute the determinant $$\begin{array}{ll}
\zz\rho\ww\rho-|\zw\rho|^2&=\frac{1}{2}v^2(\frac{1}{2}+x^2+y^2)-x^2v^2-16y^2v^2
=\frac{1}{2}(\frac{1}{2}-x^2-y^2)v^2
\end{array}$$
For small $(x,y)$ this is positive, and strictly positive if $v\ne
0$
\par We now check that
$\tilde\Phi(z,w)=P(x^2,y^2,u)+M(x^2+y^2)u^6+\frac{1}{2}v^2$ is
also plurisubharmonic. We have
$$4\ww{\tilde\Phi} = \triangle_{u,v}P+30M(x^2+y^2)u^4 + 1>
\frac{1}{2},$$ as long as $(z,w)$ is small enough. So we have an
estimate of the determinant
$$\frac{\partial^2\tilde\Phi}{\partial z\partial\bar
z}\frac{\partial^2 \tilde\Phi}{\partial w\partial \bar
w}-\left|\frac{\partial^2\tilde\Phi}{\partial z\partial\bar
w}\right|^2>\frac{1}{8}\frac{\partial^2 \tilde\Phi}{\partial
z\partial\bar z}- \left|\frac{\partial^2 \tilde\Phi}{\partial
z\partial\bar w}\right|^2.$$ \par From computations above, we know
that $\zz{\tilde\Phi}$ is bounded from below by a strictly
positive homogeneous polynomial
$$Q_5(x^2,y^2,u)=\sum a_{\kappa,\lambda,\mu}(x^2)^\kappa(y^2)^\lambda u^\mu$$
of degree $5$. We have chosen $P$ nondegenerate enough, so that
all $a_{\kappa,\lambda,\mu}$ with $\mu$ even are strictly
positive. Let us also estimate the term $|4\zw \Phi|^2$, using
(\ref{changelevi}).
$$\begin{array}{ll}|4\zw \Phi|^2
&=(\frac{\partial^2 P}{\partial u\partial
x}-(\alpha+1)x\frac{\partial^2 P}{\partial
u^2}+12Mxu^5-30M(\alpha+1)x(x^2+y^2)u^4)^2\\
&+ (-\frac{\partial^2 P}{\partial u\partial
y}+(\alpha-1)y\frac{\partial^2 P}{\partial
u^2}+12Myu^5-30M(\alpha-1)y(x^2+y^2)u^4)^2< Q_9,
\end{array}$$
where $Q_9(x^2,y^2,u)$ is some homogeneous polynomial of degree
$9$. Together, we have $$\frac{\partial^2\tilde\Phi}{\partial
z\partial\bar z}\frac{\partial^2 \tilde\Phi}{\partial w\partial
\bar w}-\left|\frac{\partial^2\tilde\Phi}{\partial z\partial\bar
w}\right|^2\ge Q_5-Q_9.$$ Since $Q_5$ is nondegenerate enough, the
polynomial $Q_5-Q_9$ is positive, as long as $(x,y,u)\ne (0,0,0)$.
Note that the sum of terms with odd powers of $u$ in $Q_9$ can be
written as $uQ_8$, where $Q_8$ is homogeneous of degree $8$, and
are thus also dominated by $Q_5$. Since $\Phi=\tilde\Phi+\rho$,
the function $\Phi$ is strictly plurisubharmonic away from the
origin. This completes the proof.
\end{proof}
Using exactly the same method as above, we can also prove the
following lemma.
\begin{lem}\label{L11} For $\alpha \le 0.44$, and any $M>0$,
there exists a homogeneous polynomial $P\in\RR[x^2,y^2,u]$ of
degree $4$ and an open neighborhood $U\subset \CC^2$ of $(0,0)$,
so that the function $$
\Phi(x,y,u)=M(x^2+y^2)u^4+P(x^2,y^2,u)+(1+x^2+y^2)v^2$$ has
properties
\begin{enumerate}\item[(a)]$\frac{\partial^2\Phi}{\partial z\partial
\bar {z}}>0,\quad (x,y,u)\in U\backslash\{(0,0)\}$ \item[(b)]
$(u,v)\cdot(\frac{\partial\Phi}{\partial
u},\frac{\partial\Phi}{\partial u})>0,\quad (z,w)\in
U\backslash\{u=v=0\}$, \item[(c)] $\Phi>0,\quad (z,w)\in
U\backslash\{u=v=0\},$ \item[(d)] $\{\Phi=0\}=\{u=v=0\}\cap U$.
\end{enumerate}
\end{lem}
This might be of some interest, since the lower the degree of the
polynomial defining the neighborhoods, the lower the vanishing of
Hessian at the origin, and thus one might expect the smaller the
irregularity of the thickness of the neighborhoods. \begin{rem}
One might try taking higher degrees of the polynomials in Lemma
\ref{L1} to push the parameter space of $\alpha$ all the way to
$1$. This would quickly make the construction very messy. We will
rather use a different construction.
\end{rem}
\par We now modify the construction above, to prove an analog of Lemma \ref{L1} for
$\alpha>\frac{1}{2}$. We want to use a construction, similar to
the one in Lemma \ref{L1}, but having fewer inequalities and
equalities to satisfy. We will achieve that, by adding a function
that has a sufficiently large Levi form, and vanishes on $u=0$.
\begin{lem} \label{L2} Let $\frac{1}{2}<\alpha<1$. Then there exists a
positive integer $m$, such that for every $\delta>0$ the function
$\Phi_m(z,w):=y^{2m}u^2+\delta v^2$ is plurisubharmonic in a small
neighborhood $U$ of the origin, and strictly plurisubharmonic in
$U\backslash\{y=0\}$.
\end{lem}
\begin{proof} Let $\alpha<1$ be fixed. We check the condition for
the function $\phi_m=y^{2m}u^2$ to have the property
$\zz{\phi_m}>0$ for small $(z,w)$ and $y\ne 0$.
\begin{equation}\begin{array}{ll} \label{temp1}
4\zz{\phi_m}=&2m(2m-1)y^{2m-2}u^2-4((2m+1)\alpha-2m)y^{2m}u+2(\alpha-1)^2y^{2m+2}\\
            &+2(\alpha+1)^2y^{2m}x^2.
\end{array}\end{equation} We get the desired positivity if and only if
\begin{equation}\nonumber
2m(2m-1)y^{m-2}u^2-4((2m+1)\alpha-2m)y^mu+2(\alpha-1)^2y^{2m+2}
\end{equation}
is positive for $y\ne 0$. Using Lemma \ref{pos}, the condition is
equivalent to $$\begin{array}{ll}
((2m+1)\alpha-2m)^2&<m(2m-1)(\alpha-1)^2
\\&\Leftrightarrow
\frac{2m-\sqrt{(2m-1)m}}{2m+1-\sqrt{(2m-1)m}}<\alpha<\frac{2m+\sqrt{(2m-1)m}}{2m+1+\sqrt{(2m-1)m}}.
\end{array}$$
Let $a_m:=\frac{2m-\sqrt{(2m-1)m}}{2m+1-\sqrt{(2m-1)m}}$ and
$b_m:=\frac{2m+\sqrt{(2m-1)m}}{2m+1+\sqrt{(2m-1)m}}$. Thus
(\ref{temp1}) holds for some $m$ and $a_0<\alpha<\lim_{m\to\infty
}b_m=1$, if we show
\begin{enumerate}
\item $a_m,b_m$ are both increasing sequences, \item
$a_{m+1}<b_m$.
\end{enumerate}
We first show that $a_m$ is an increasing sequence.
$$\begin{array}{ll}&a_m<a_{m+1}\\ &\Longleftrightarrow
\left({2m+2-\sqrt{(2m+1)(m+1)}}\right)\left({2m+1-\sqrt{(2m-1)m}}\right)>\\
&\left({2m+3-\sqrt{(2m+1)(m+1)}}\right)\left({2m-\sqrt{(2m-1)m}}\right)\\
&\Longleftrightarrow
2>\sqrt{(2m+1)(m+1)}-\sqrt{(2m-1)m}\Leftrightarrow 16m^2+8m-9>0.
\end{array}$$
Showing  $b_m$ is increasing is exactly the same, but it also
follows from
$$\begin{array}{ll}&a_{m+1}<b_m\\ &\Longleftrightarrow
\left({2m+\sqrt{(2m-1)m}}\right)\left({2m+3-\sqrt{(2m+1)(m+1)}}\right)>\\
&\left({2m+1+\sqrt{(2m-1)m}}\right)\left({2m+2-\sqrt{(2m+1)(m+1)}}\right)\\
&\Longleftrightarrow 2<\sqrt{(2m+1)(m+1)}+\sqrt{(2m-1)m}.
\end{array}$$
This obviously holds for $m\ge 1$. Since $a_1=\frac{1}{2}$, we
have (\ref{temp1}) for $\frac{1}{2}<\alpha<1$.
\par We now show that $\Phi_m$ is plurisubharmonic in a small neighborhood of $(0,0)$,
as long as $m$ is chosen so that $a_m<\alpha<b_m$. We already have
shown that $\zz{\Phi_m}=\zz{\phi_m}$ is nonnegative. For $y$
small, we have
$$\zz{\Phi_m}\ww{\Phi_m}-\left|\zw{\Phi_m}\right|^2>\frac{\delta}{2}\zz{\phi_m}-\left|\zw{\phi_m}\right|^2.
$$
Similar than in Lemma \ref{L1}, monomials in
$\left|\zw{\Phi_m}\right|^2$ are dominated by monomials in
$\frac{\delta}{2}\zz{\phi_m}$. This completes the proof of this
lemma.
\end{proof}
We now combine our first approach with Lemma \ref{L2} above, to
get desired functions for all quadratic hyperbolic points. Let us
first prove this elementary lemma.
\begin{lem}\label{pos2} The function
$$q(x,y,u)=u^{2k}+a|x|^\gamma|y|^\delta u^l+b|x|^{\gamma_1}|y|^{\delta_1}$$
is strictly positive for small $(x,y,u)\ne (0,0,0)$, as long as
$b>0$, $\gamma_1<\frac{2k\gamma}{2k-l}$,
$\delta_1<\frac{2k\delta}{2k-l}$.
\end{lem}
\begin{proof} For fixed $(x,y)$, the minimum of our
function is achieved at the critical point:
$$\frac{\partial q}{\partial
u}=u^{l-1}(2ku^{2k-l}+al|x|^\gamma|y|^\delta)=0.$$ If the minimum
is obtained at $u=0$, we are done. If not, then at the minimum
$$u=o(|x|^\frac{\gamma}{2k-l}|y|^\frac{2k}{2k-l}),$$  and so
$$\min_uq(x,y,u)=o(|x|^\frac{2k\gamma}{2k-l}|y|^\frac{2k\delta}{2k-l})+b|x|^{\gamma_1}|y|^{\delta_1}.$$
The value is strictly positive for $x,y$ nonzero and small.
\end{proof}

\begin{lem}\label{L3}  Let $1>\alpha > \frac{1}{2}$. Then there
exist positive integers $n,m,k$, a homogeneous polynomial
$P\in\RR[x^2,y^2,u]$ of degree $2n$, such that for every $M>0$,
the function
\begin{eqnarray}\nonumber
\Phi(z,w)=P(x^2,y^2,u)+y^{2m}u^2+x^{2k}u^2+M(x^2+y^2)u^{2n}+(1+x^2+y^2)v^2
\end{eqnarray}
has the following properties on a small neighborhood $U$ of
$(0,0)\in\CC^2$
\begin{enumerate} \item[(a)] $\Phi$ is plurisubharmonic in $U$ and
strictly plurisubharmonic in $U\backslash\{(0,0)\}$, \item[(b)]
$(u,v)\cdot(\frac{\partial\Phi}{\partial
u},\frac{\partial\Phi}{\partial v})>0,\quad (z,w)\in
U\backslash\{u=v=0\}$, \item[(c)] $\Phi>0,\quad (z,w)\in
U\backslash\{u=v=0\}$, \item[(d)] $\{\Phi=0\}=\{u=v=0\}\cap U$.
\end{enumerate}
\end{lem}
\begin{proof} The proof is very similar to the proof
of Lemma \ref{L1}. Let us choose $\frac{1}{2}<\alpha<1$. We set
$m$ as in Lemma \ref{L2}, so that $\zz{(y^{2m}u^2)}$ is strictly
positive for $y\ne 0$. Let $n$ be any positive integer with
$2n>m+2$ and $k$ any integer with $k+2>2n$. Let
$P\in\RR[x^2,y^2,u]$ be a homogeneous polynomial of degree $2n$ of
the form
$$u^{2n}+(ax^2+cy^2)u^{2n-1}+(Ax^4+Bx^2y^2+Cy^4)u^{2n-2},$$
where $a,c,A,B,C$ are going to be appropriately chosen. Using
Lemmas \ref{pos} and \ref{pos2}, we get (a) and (b) to hold for
$\Phi$ as long as
\begin{enumerate}
\item[1.] $B>0$ \item[2.] $8(2n-2)nA>(2n-1)^2a^2$.
\end{enumerate}
Calculating $\zz{P}$ we get
$$\begin{array}{ll}\zz{P}  &=[2a+2b-4n\alpha]u^{2n-1}\\
&+[(12A+2B+2n(2n-1)(\alpha+1)^2-2(2n-1)a(3\alpha+2))x^2\\
&+(12C+2B+2n(2n-1)(\alpha-1)^2-2(2n-1)c(3\alpha-2))y^2]u^{2n-2}\\
&+[((2n-1)(2n-2)a(\alpha+1)^2-2(2n-2)(5\alpha+4)A)x^4\\
&+((2n-1)(2n-2)a(\alpha-1)^2+(2n-1)(2n-2)c(\alpha-1)^2-10(2n-2)B)x^2y^2\\
&+((2n-1)(2n-2)c(\alpha-1)^2-2(2n-2)(5\alpha-4)C)y^4]u^{2n-3}\\
&+[(Ax^4+Bx^2y^2+Cy^4)((\alpha+1)^2+(\alpha-1)^2)]u^{2n-4}
\end{array}$$
We proceed by setting some of the terms equal to $0$ and some of
the terms to be positive, but we will see later, that we do not
have to worry about what happens to terms expressed purely by $y$
and $u$. So we just want the following equalities and inequalities
to hold
\begin{enumerate}
\item[3.] $a+c=2n\alpha$ \item[4.]
$(2n-1)a(\alpha+1)^2=2(5\alpha+4)A$ \item[5.]
$(2n-1)a(\alpha-1)^2+(2n-1)c(\alpha+1)^2=10(2n-2)B$ \item[6.]
$12A+2B+2n(2n-1)(\alpha+1)^2>2(2n-1)a(3\alpha+2)$
\end{enumerate}
As long as $a$ is chosen positive and small enough, and we set
$c=2n\alpha-a$, all equalities and inequalities $1.-6.$ are
satisfied. Choosing also $C$ large enough, we have that
$$\zz P=(\tilde A x^2 +\tilde C y^2)u^{2n-2}+\tilde D y^4 u^{2n-3}+
(\tilde Ex^4+\tilde F x^2y^2+\tilde G y^4),$$ where $\tilde
A,\tilde B,\tilde E ,\tilde F ,\tilde G>0$. Putting it all
together, we get $$\begin{array}{ll}&\zz
\Phi=M\zz{((x^2+y^2)u^2)}+\zz{(y^{2m}u^2)}+\zz{(x^{2k}u^2)}+\zz{P}+\zz{((x^2+y^2)v^2}\ge\\
&4Mu^{2n}-4nM((3\alpha+2)x^2+(3\alpha-2)x^2)u^{2n-1}+(\tilde A x^2
+\tilde C y^2)u^{2n-2}\\&+\tilde D
y^4u^{2n-3}-4((2k+1)\alpha+2k)x^{2k}u+2(\alpha+1)^2x^{2k+2}+\zz{(y^{2m}u^2)}\\&=
[Mu^{2n}-4nM((3\alpha+2)x^2+(3\alpha-2)x^2)u^{2n-1}+(\tilde A
x^2+\tilde C y^2)u^{2n-2}]\\
&+[Mu^{2n}-4((2k+1)\alpha+2k)x^{2k}u+2(\alpha+1)^2x^{2k+2}]\\
&+[Mu^{2n}+\tilde D y^4u^{2n-3}+\zz{(y^{2m}u^2)}]+Mu^{2n}.
\end{array}$$
Each of the summands in the last equality are positive, by using
Lemma \ref{pos}.
\par To prove that the determinant
$$\zz\Phi\ww\Phi-\left|\zw\Phi\right|^2$$
is positive, we proceed exactly the same way as in Lemma \ref{L1},
by comparing the orders of monomials. This proves the lemma.
\end{proof}
Putting things together, we have proven the following proposition.
\begin{prop} Let $S$ be a real surface in $\CC^2$ given by the equation
$$w=\frac{1}{2}\alpha z\bar z+\frac{1}{4}(z^2+{\bar z}^2),\qquad 0\le\alpha<1.$$
Then there exists a $\C^\infty$ function $\Phi$, defined in a
small neighborhood $U$ of $(0,0)\in \CC^2$ with properties
\begin{enumerate}
\item[(a)] $\Phi$ is plurisubharmonic in $U$ and strictly
plurisubharmonic in $U\backslash\{(0,0)\}$, \item[(b)]$\Phi\ge
0$,\item[(c)] $\{\Phi=0\}=\{\nabla \Phi =0\}=S\cap U$.
\end{enumerate}
\end{prop}
\subsection{The flat case.}
\begin{lem} Let $\pi\!:S\hookrightarrow X$ be a real surface
imbedded in a complex surface $X$. Let $p$ on $S$ be a hyperbolic
flat complex point on $S$. Then there exists a neighborhood $V$ of
$p$ in $X$ and a smooth function $\phi\!:V\to \RR$ with properties
\begin{enumerate}
\item[(a)] $S\cap V=\{\phi=0\}=\{\nabla\phi\}=0$ \item[(b)]
$\phi\ge 0$ \item[(c)] $\phi$ is strictly plurisubharmonic in
$V\backslash \{p\}$.
\end{enumerate}
\end{lem}
\begin{proof} Let $U$ be a neighborhood of $p$ and $(z,w)$ be holomorphic coordinates
on $U$ in which $S$ is written as in (\ref{formatocke}), where
$0\le\alpha<1$ and $\tau=o(|z|^3)$. We use the non-holomorphic
coordinates (\ref{nonholomorphic}) to compute the Levi form. Let
first $\alpha\le \frac{1}{2}$, and let $P\in\RR[x^2,y^2,u]$ be the
homogeneous polynomial, constructed in Lemma \ref{L1}. We take
$$\Phi(z,w)=P(x^2,y^2,u)+M(x^2+y^2)u^6+v^2,$$ where
$M$ is to be chosen later. Similar calculations as in Lemma
\ref{L1} show that
$$4\zz \Phi=4Mu^6+q_1u^5+q_2u^4+q_5u^3+q_6u^2+q_9u+q_{10}$$
where $q_2,q_6,q_{10}\in \RR[x,y]$ are positive polynomials
respectively of degrees $2,6,10$, vanishing only at the origin,
and $q_1,q_5,q_9\in \RR[x,y]$ are some polynomials respectively of
degrees $1,5,9$. These last terms are the error terms, coming from
$$((o(2)\frac{\partial }{\partial x}+o(2)\frac{\partial }{\partial
y}+o(1))\frac{\partial }{\partial u}+o(3)\frac{\partial^2
}{\partial u^2})\Phi.$$ Using Lemma \ref{pos}, $\zz \Phi$ is
positive away from the origin, as long as $q_1^2<16Mq_2$. Choosing
$M$ large enough, this is satisfied. Since we have that $4\ww
\Phi=2+o(1)$, the expression $\zz\Phi\ww\Phi$ approximately equals
$\frac{1}{2}\zz\Phi$ near the origin. Comparing degrees of
monomials as we did in Lemma \ref{L1}, we see that
$$\zz\Phi\ww\Phi-\left|\zw \Phi\right|^2>0$$
away from the origin.
\par We are left with $\alpha>\frac{1}{2}$. In this case, we use
the function $\Phi$ constructed in Lemma \ref{L3}. The rest
follows exactly the same way.
\end{proof}
\section{Stein neighborhoods of real surface}
Using the results from the previous section, we are ready to prove
Theorem \ref{maintheorem}. First we need this classical result.
\begin{prop} Let $S^n$ be a $\C^\infty$ smooth totally real regularly embedded submanifold  in a complex
manifold $X^n$. Let $<\cdot,\cdot>$ be any $\C^\infty$ smooth
Riemannian metric on $X$. The function $\Phi\!:X\to \RR$, defined
as
$$\Phi(z)=dist^2(z,S)$$
is $\C^\infty$ smooth and strictly plurisubharmonic in an open
neighborhood of $S$ in $M$.
\end{prop}
\begin{proof}Let $p\in S$ be any point, and let $Z$ be a
smooth vector field, defined in a neighborhood of $p$ in $X$, with
$Z(p)\ne 0$. We need to show $d\dc\Phi$ is a positive form in a
neighborhood of $S$. Let $J$ be the complex structure on $X$.
$$\begin{array}{ll}
d\dc \Phi(p)(Z,JZ)&=(Z\dc\Phi(JZ))(p)-((JZ)\dc\Phi(Z))(p)-(\dc\Phi([Z,JZ]))(p)\\
&=(Zd\Phi(Z))(p)+((JZ)d\Phi(JZ))(p)+(d\Phi(J[Z,JZ]))(p)\\&=<\Hess\Phi(p)
Z,Z>_p+<\Hess\Phi(JZ),(JZ)>_p,\end{array}$$ where $\Hess\Phi$ is
the real Hessian form. Since $S$ is maximally real, either $Z$ or
$JZ$ must have a normal component. Since $S$ is regularly
embedded, $\Hess\Phi$ is a non-degenerate positive form on the
normal bundle of $S$ in $X$. This completes the proof.
\end{proof}
\begin{thm} \label{T1} Let $S\hookrightarrow X$ be a compact real surface,
$\C^\infty$ embedded into a complex surface $X$ and having only
flat hyperbolic complex points $\{p_1,\ldots,p_k\}$. Then there
exists a $\C^\infty$ function $\psi$, defined in a neighborhood
$U$ of $S$ in $X$, such that
\begin{enumerate}
\item[(a)] $S=\{\psi=0\}=\{\nabla\psi=0\}$, \item[(b)] $\psi$ is
strictly plurisubharmonic on $U\backslash\{p_1,\ldots,p_k\}$.
\end{enumerate}
Sublevel sets $\Omega_\epsilon=\{\psi<\epsilon\}$ define a
regular, strictly pseudoconvex Stein neighborhood basis of $S$ in
$M$. \end{thm}
\begin{proof}
Let $\tilde S=S\backslash\{p_1,\ldots,p_k\}$. For every $1\le j\le
k$ let $\Phi_j\!:U_j\to \RR$ be the plurisubharmonic function,
constructed in the previous section and defined in a small
neighborhood $U_j$ of $p_j$. All $U_j$ are assumed to be pairwise
disjoint. If $U_j$ is taken small enough, $\Phi_j$ is strictly
plurisubharmonic in $U_j\backslash\{p_j\}$. Let $<\cdot,\cdot>$ be
any Riemannian metric on $X$ and let $\Phi_0(z)=dist^2(z,S)$. Then
$\Phi$ is strictly plurisubharmonic in a neighborhood $U_0$ of
$\tilde S$ in $X$. Let $V=\bigcup_{0\le j\le k}U_j$ and let
$\pi\!:V\to S$ be the map defined as $\pi(z)=p$,
$dist(z,p)=dist(z,S)$. Provided the neighborhoods are chosen small
enough, the map $\pi$ is well defined and $\C^\infty$. Let
furthermore $\{\rho_j\}_{0\le j\le k}$ be a $\C^\infty$ partition
of unity for $\{U_j\cap S\}_{0\le j\le k}$. We define
$$\Phi(z):=\sum_{j=0}^k\rho_j(\pi(z))\Phi_j(z).$$ For every $p\in S$, we have
$$d\dc\Phi(p)=\sum_0^k\rho_j(p) d\dc \Phi_j(p).$$ This
expression is strictly positive away from the points
$\{p_1,\ldots,p_j\}$. We also have that $d\dc\Phi=d\dc\Phi_j$ near
$p_j$. After possibly shrinking $V$, $\Phi$ is plurisubharmonic in
a neighborhood $U$ of $S$ in $X$ and strictly plurisubharmonic in
$U\backslash \{p_1,\ldots,p_k\}$. Since
$\nabla\Phi=\nolinebreak\nabla\Phi_j$ near $p_j$, we also have
$\nabla\Phi$ nonvanishing near $S$. \par What is left is to show
that the neighborhoods are indeed Stein. By a result of Grauert
\cite{Gra1}, strictly pseudoconvex domains in a complex manifold
are Stein, if and only if they contain no compact complex analytic
sets of positive dimension. The restrictions of strictly
plurisubharmonic functions to analytic sets are again strictly
plurisubharmonic. Since compact analytic sets do not have any
nonconstant plurisubharmonic functions and the defining function
$\Phi$ is strictly pseudoconvex everywhere, but at finitely many
points, there can not be any compact positive dimensional analytic
sets in our neighborhood. This completes the proof.
\end{proof}
The above theorem, together with Lemma \ref{mama lemma}, proves
Theorem \ref{maintheorem}. \par It would be nice to know if such
functions can be constructed without the assumption of flatness of
hyperbolic points. For now, we satisfy ourself with the next
result.
\begin{cor}\label{C1} Let $\pi\!:S\hookrightarrow X$ be any
generically embedded real compact surface without elliptic points
in a complex surface $X$. Then there exists an embedding\hfill \\
$\pi'\!:S\hookrightarrow X$, $\C^{2}$ close to $\pi$ and isotopic
to $\pi$, such that $\pi'(S)$ has a regular basis of Stein
neighborhoods.
\end{cor}
\begin{proof} In local coordinates $(z,w)$ near a hyperbolic complex point
$p$, the surface $S$ can be written as
$$\RE w=\frac{1}{2}\alpha z\bar z+\frac{1}{4}(z^2+{\bar z}^2)+\tau_1(z),\qquad \IM w=\tau_2(z)$$
where $\tau_1,\tau_2=o(|z|^3)$ and real. Let $\rho\!:S\to [0,1]$
be a smooth function with $\rho\equiv 1$ near $p$ and $\rho\equiv
0$ outside a small neighborhood of $p$ . Then
$\pi_t:=\pi(t)-it\rho\tau_2$ defines an isotopy of $\pi=\pi_0$ to
$\pi'=\pi_1$ with a flat complex point at $p$. Repeating this
process for every complex point, and using Theorem \ref{T1} we get
the required result.
\end{proof}
\section{Application to unions of totally real planes in $\CC^2$}
Using results from the previous sections, we construct Stein
neighborhood basis for certain unions of two totally real planes
$L_1,L_2\subset \CC^2$, with $L_1\cap L_2=\{0\}$. Every such union
is linearly holomorphically equivalent to $M(B):=\RR^2\cup A(B)$,
where $A(B)$ is the real span of the columns of the matrix $B+iI$,
with $B$ real and $(B-iI)$ invertible. Furthermore, $B$ is
determined only up to real conjugancy. For the proofs of these
simple statements, see \cite{Wei}.
\par The next two lemmas show the connection
between certain unions of totally real planes and complex points
on real surfaces. We only prove the first one, since it is the one
we use later.
\begin{lem}\label{hyperbolic} Let $B=\begin{pmatrix}
  \mu & 0 \\
  0 & -\mu \\
\end{pmatrix},$ with $\mu\ge 0$.  Let
$$\Psi(z,w)=\left(ie^{-i\frac{\theta}{2}}z+ie^{i\frac{\theta}{2}}w,-(\sin^2{\theta})zw\right),$$
where $\theta=\arcsin\sqrt{\frac{1}{1 +\mu^2}}$. Then
$$\Psi(M(B))=\{(z,w),w=\frac{1}{2}\alpha z\bar
z+\frac{1}{4}z^2+\frac{1}{4}{\bar z}^2\},$$ with $\alpha=\cos
\theta$. \phantom{la}
\end{lem}
\begin{proof} The set $\Psi(M(B))$ can be parameterized by
$$(ie^{\frac{-i\theta}{2}}(i+\mu)x+
ie^{\frac{i\theta}{2}}(i-\mu)y,-(\sin^2
\theta)(i+\mu)(i-\mu)xy).$$ Let us use this to compute
$\frac{1}{2}\alpha z\bar z+\frac{1}{4}z^2+\frac{1}{4}{\bar z}^2$
on $\Psi(M(B))$.
$$\begin{array}{ll}
\frac{1}{2}\alpha z\bar z&+\frac{1}{4}z^2+\frac{1}{4}{\bar z}^2=\\
&\frac{1}{2}\alpha |ie^{\frac{-i\theta}{2}}(i+\mu)x+
ie^{\frac{i\theta}{2}}(i-\mu)y|^2+\RE{\frac{1}{2}(ie^{\frac{-\theta}{2}}(i+\mu)x+
ie^{\frac{\theta}{2}}(i-\mu)y)^2}=\\
&[\frac{1}{2}\alpha(\mu^2+1)-\frac{1}{4}(\mu+i)^2e^{-i\theta}-\frac{1}{4}(\mu-i)^2e^{i\theta}](x^2+y^2)\\
&+[\frac{1}{2}\alpha(-(\mu+i)^2e^{-i\theta}-(\mu-i)^2e^{i\theta})+(\mu^2+1)]xy\\
&=xy\end{array}$$ In the above, we used that
$\alpha=\frac{1}{\sqrt{\mu^2+1}}$. On the other hand, the second
coordinate on $\Psi(M(B))$ equals
$$-(\sin^2 \theta)(i+\mu)(i-\mu)xy=-\frac{1}{1+\mu^2}(-1-\mu^2)xy=xy.$$
So we have concluded that
$$\frac{1}{2}\alpha z\bar
z+\frac{1}{4}z^2+\frac{1}{4}{\bar z}^2=w$$ on $\Psi(M(B))$. This
concludes the proof.
\end{proof}
\begin{lem}\label{elliptic} Let $B=\begin{pmatrix}
  0 & \mu \\
  -\mu & 0 \\
\end{pmatrix}$ with $\mu>1$. Let
$$\Psi(z,w)=\left(ie^{-i\frac{\theta}{2}}z+ie^{i\frac{\theta}{2}}w,
\frac{1}{2}(\tan\theta \sin\theta)(z^2+w^2)\right),$$ where
$\theta=\arcsin{\frac{1}{\mu}}$. Then
$$\Psi(M(B))=\{(z,w),w=\frac{1}{2}\alpha z\bar
z+\frac{1}{4}z^2+\frac{1}{4}(\bar z)^2\},$$ with
$\alpha=\frac{1}{\cos \theta}$.
\end{lem}
\begin{rem}
The maps constructed above were, in a slightly different context,
found by Burns (personal communication). Lemma \ref{elliptic} can
be used to pull-back Bishop discs, \cite{Bis}, from a neighborhood
of an elliptic complex point. This gives us analytic annuli with
boundaries in $M(B)$ with $\trace B=0$ and $\det B>1$, shrinking
towards the origin. Weinstock \cite{Wei} showed that this is the
only non-polynomially convex case among the unions $M(B)$. One
would thus expect to be able to find a regular Stein neighborhood
basis for all other unions $M(B)$. Unfortunately, we are at the
point only able to show this for a smaller class of unions of
totally real planes.
\end{rem}
\begin{prop} Let $B$ be a real $2\times 2$ matrix, diagonalizable over
$\RR$ and with the property that $\trace B =0$. Then the union
$M(B)$ has a regular Stein neighborhood basis.
\end{prop}
\begin{proof} We can assume $B=\begin{pmatrix}
  \mu & 0 \\
  0 & -\mu \\
\end{pmatrix}$. By Lemma \ref{hyperbolic}, the map
$$\Psi(z,w)=\left(ie^{-i\frac{\theta}{2}}z+ie^{i\frac{\theta}{2}}w,-(\sin^2{\theta})zw\right)$$
with $\theta=\arcsin\sqrt{\frac{1}{1+\mu^2}}$ maps the union
$M(B)$ of totally real planes onto the surface
$$S_\alpha=\{(z,w),w=\frac{1}{2}\alpha z\bar
z+\frac{1}{4}z^2+\frac{1}{4}{\bar z}^2\}.$$ Let $\Phi\!:U\to \RR$
be the map, constructed in Lemma \ref{L1} or Lemma \ref{L3},
depending on the size of $\alpha=\cos\theta$. Let
$\tilde\Phi=\Phi\circ\Psi$. The small sublevel sets of
$\tilde\Phi$ are pseudoconvex. We only need to check that
$\nabla\tilde\Phi(z,w)\ne 0$, for $(z,w)\ne M(B)$. Since
$\nabla\Phi$ is nonzero away from the surface $S_\alpha$, this
happen if $\nabla\Phi\notin \ker (D\Psi)^T$. We have
$$D\Psi=\begin{pmatrix}
 ie^{-i\frac{\theta}{2}}  & ie^{i\frac{\theta}{2}}\\
  -(\sin^2{\theta})w & -(\sin^2{\theta})z \\
\end{pmatrix}.$$ $D\Psi$ is nondegenerate outside of
$z=e^{i\theta}w$ and $\ker
D\Psi(e^{i\theta}w,w)=\CC\{(e^{i\frac{\theta}{2}},-e^{-i\frac{\theta}{2}})\}$.
Let us assume that, at some point $(z_0,w_0)\ne (0,0)$ in the
image of the critical set, $\nabla\Phi\in\ker (D\Psi)^T$. We know
from constructions of $\Phi$, that $$\Phi+N|z|^2(Re w
-\frac{1}{2}\alpha z\bar z-\frac{1}{4}z^2-\frac{1}{4}{\bar
z}^2)^{2n}$$ is also plurisubharmonic. Here $2n$ is the degree of
the homogeneous polynomial in the definition of $\Phi$. So by
possibly substituting $\Phi$ by $$\Phi +N|z|^2(Re w
-\frac{1}{2}\alpha z\bar z-\frac{1}{4}z^2-\frac{1}{4}{\bar
z}^2)^{2n}$$ for an appropriate $N$, ensures us to have no
critical points of $\tilde\Phi$ away from $M(B)$. This completes
the proof.
\end{proof}
\subsection*{Acknowledgment} I would like to thank John Erik Forn\ae ss,
for many stimulating discussions on the topic.

\end{document}